\documentclass{style}
\usepackage{cite}

\DeclareMathOperator{\vol}{vol}
\DeclareMathOperator{\sol}{sol}
\DeclareMathOperator{\tsol}{tsol}
\DeclareMathOperator{\B}{B}
\DeclareMathOperator{\density}{\delta}
\DeclareMathOperator{\G}{G}
\DeclareMathOperator{\voronoi}{\Omega}
\DeclareMathOperator{\funcL}{L}
\DeclareMathOperator{\h}{h}
\DeclareMathOperator{\aff}{aff}
\DeclareMathOperator{\dimaff}{dim aff}
\DeclareMathOperator{\card}{card}
\DeclareMathOperator{\conv}{conv}
\DeclareMathOperator{\dih}{dih}
\DeclareMathOperator{\cell}{cell}
\DeclareMathOperator{\M}{M}
\DeclareMathOperator{\EC}{EC}
\DeclareMathOperator{\wt}{wt}
\DeclareMathOperator{\CL}{CL}

\begin{document}
\title{On a Detail in Hales's "Dense Sphere Packings: A Blueprint for Formal Proofs" }

\author{Nadja Scharf%
  \thanks{\texttt{nadja.scharf@fu-berlin.de}}}
\affil{Institute of Computer Science, Freie Universität Berlin, Takustr.~9, 14195~Berlin, Germany}
	
	\maketitle
	\begin{abstract}
	In \cite{Blueprint} Hales proves that for every packing of unit spheres, the density in a ball of radius $r$ is at most $\pi/\sqrt{18}+c/r$ for some constant $c$. When $r$ tends to infinity, this gives a proof to the famous Kepler conjecture. As formulated in \cite{Blueprint}, $c$ depends on the packing. We follow the proofs in \cite{Blueprint} to calculate a constant $c'$ independent of the sphere packing that exists as mentioned in \cite{Hales}.
	\end{abstract}
\section{Introduction}
\label{sec:Introduction}	
In \cite{Blueprint} Hales proves that for every packing of infinitely many unit spheres into three dimensional space, there exists a constant $c$ such that the \emph{density} inside a ball of radius $r$ is upper bounded by $\pi/\sqrt{18}+c/r$. Here, the density is defined as the volume of the intersection of the packed unit spheres with the container sphere of radius $r$ divided by the volume of the container sphere. The famous Kepler conjecture states that the density tends to $\pi/\sqrt{18}$ when $r$ tends to infinity which is implied by the density bound shown by Hales.

To make, for example, statements about bounds for finite packings inside a container, we need $c$ to be independent of the packing, i.e., we want a statement of the form: There exists a constant $c'$ such that for every infinite sphere packing the density inside a ball of radius $r$ is upper bounded by $\pi/\sqrt{18}+c'/r$. If it holds for an infinite packing, then this statement would also hold for a finite packing.

First, we give some definitions that are necessary to understand the crucial lemmas from \cite{Blueprint}. Then, we follow the proofs of the lemmas and calculate a constant $c'$ independent of the packing. Instead of using $\mathcal{O}$-Notation as in \cite{Blueprint}, we give more detailed calculations to be able to give an actual value for $c'$. When we cite lemmas or definitions from \cite{Blueprint} (with occasional slight modifications) we give the corresponding number of the lemma or definition in \cite{Blueprint} parenthesized. We try to give definitions as closely as possible to the point where we use them.

\section{Main Lemmas}
Let us denote by $\vol$ the Lebesgue measure on Euclidean space $\mathbb{R}^3$.

In the sequel, let $V \subset \mathbb{R}^3$ be a point set that induces a packing of infinitely many unit spheres, i.e., the points in $V$ are the centers of the spheres and thus have pairwise distance at least 2. Let $\B(\mathbf{p},r)$ be the open sphere centered at point $\mathbf{p}$ with radius $r$. Let $V(\mathbf{p},r) = V \cap \B(\mathbf{p},r)$. The \emph{density} $\density (V,\mathbf{p},r)$ inside a container sphere centered at point $\mathbf{p}$ with radius $r$ is defined by 
\begin{equation*}
\density(V,\mathbf{p},r)= 
\frac{\vol \left( \B(\mathbf{p},r)\cap \bigcup_{\mathbf{v}\in V} \B(\mathbf{v},1) \right) }
{\vol\left( \B(\mathbf{p},r)\right) }.
\end{equation*}
In \cite{Blueprint}, it is shown that the face-centered cubic (FCC) packing of unit spheres has optimal density when $r$ tends to infinity. In the FCC-packing, the spheres are arranged in layers. In each layer, the spheres are arranged with there centers on a hexagonal grid where neighboring grid points have distance 2. The vertices of the second layer lie above centers of the triangles in the first layer. In the third layer, the vertices lie above centers of triangles in the first and in the second layer. These three layers are then repeated infinitely often. The layers are pushed together until the spheres touch. See the following picture for illustration.
\begin{figure}[H]
\centering
\includegraphics[width = 0.3\textwidth]{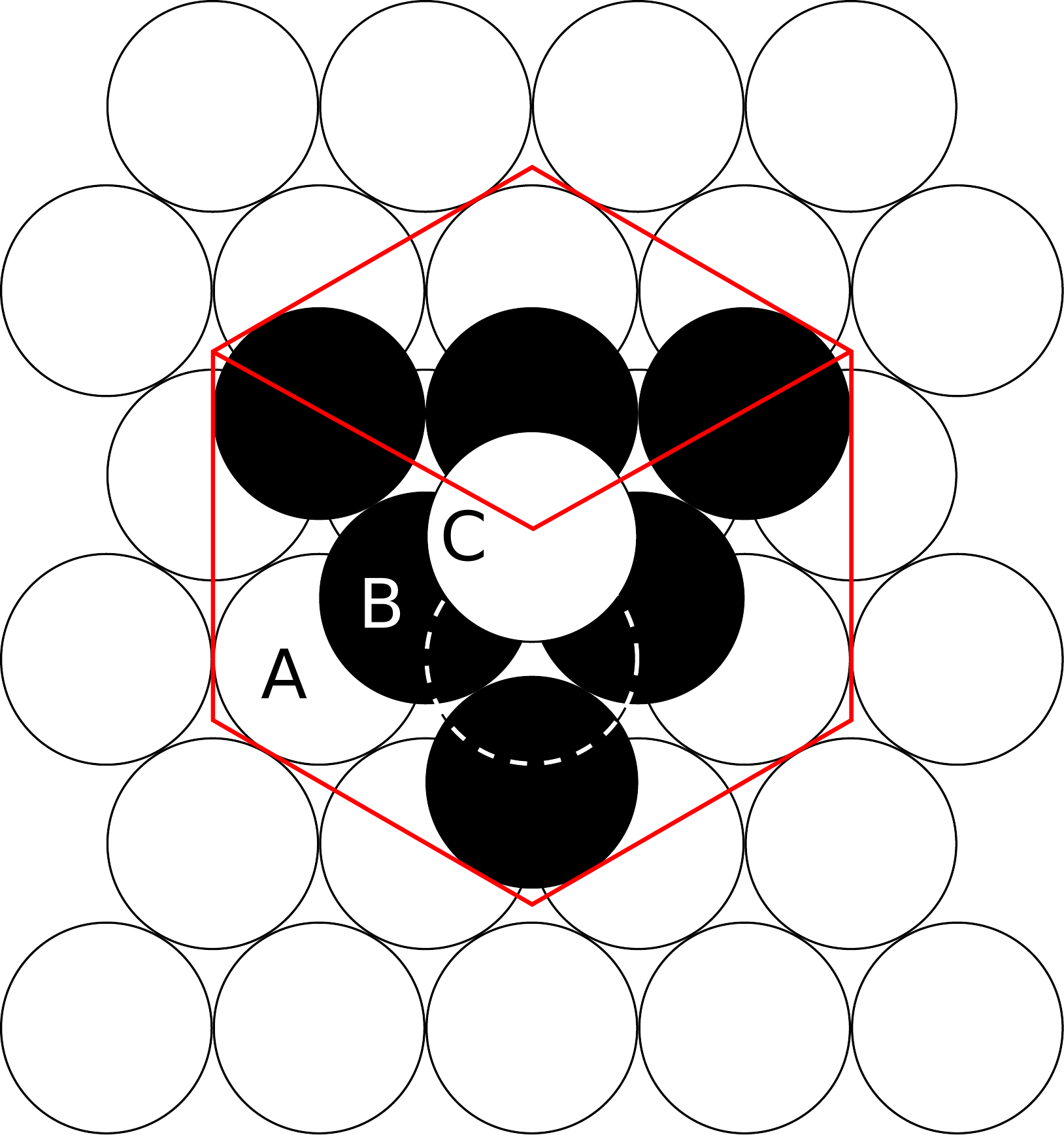}
\caption{Three layers of the FCC-packing}
\end{figure}
Consider the Voronoi diagram of the circle centers in a FCC-packing. The Voronoi cells are dodecahedra with volume $4\sqrt{2}$ which yields a density of the packing of $\pi/\sqrt{18}$.

Let $\voronoi(V,\mathbf{v})$ denote the cell of $\mathbf{v}$ in the Voronoi diagram of $V$ and, more generally, let $\voronoi(V,\underline{\mathbf{u}})$ be the intersection of the Voronoi cells for points in the list $\underline{\mathbf{u}}=[\mathbf{u}_0;\mathbf{u}_1, \dots]$.

\begin{newdef}[Definition 6.11 ]
A function $\G\colon V\to \mathbb{R}$ on a set $V \subset \mathbb{R}^3$ is \emph{negligible} if there is a constant $c_1$ such that for all $r \geq 1$, 
\begin{equation*}
\sum_{\mathbf{v} \in \B(\mathbf{0},r) \cap V} \G(\mathbf{v}) \leq c_1 r^2.
\end{equation*}
A function $\G\colon V \to \mathbb{R}$ is \emph{FCC-compatible} if for all $\mathbf{v} \in V$
\begin{equation*}
4\sqrt{2} \leq \vol(\voronoi(V,\mathbf{v}))+\G({\mathbf{v}}).
\end{equation*}
\label{def:FCCnegl}
\end{newdef}
\emph{FCC-compatible} means that the volume of every Voronoi cell is close to the volume of those in  the FCC-packing. Then, \emph{negligible} means, that the error is small.

Now we are ready to look at a lemma in \cite{Blueprint} from which the Kepler conjecture follows under certain assumption.
\begin{lem}[Lemma 6.13 from\cite{Blueprint}] 
If there exists a negligible FCC-compatible function $\G\colon V \to \mathbb{R}$ for a saturated packing $V$, then there exists a constant $c = c(V)$ such that for all $r \geq 1$ 
\begin{equation*}
\density(V,\mathbf{0},r) \leq \frac{\pi}{\sqrt{18}}+\frac{c}{r}.
\end{equation*}
\label{lem:6.13}
\end{lem}

Since this lemma is not sufficient to show the Kepler conjecture, two more lemmas are required. Therefore, we need two more definitions.

\begin{newdef}[Definition 6.34 ]
If $\underline{\mathbf{u}}=\left[\mathbf{u}_0; \dots ; \mathbf{u}_k\right]$ is a list of points in $\mathbb{R}^n$, then let $\h(\underline{\mathbf{u}})$ be the circumradius of its point set $\left\lbrace \mathbf{u}_0; \dots ; \mathbf{u}_k \right\rbrace$.
\end{newdef}

\begin{newdef}[Definition 6.88 ]
Set
\begin{equation*}
h_0 = 1.26.
\end{equation*}
Let $\funcL\colon\mathbb{R}\to \mathbb{R}$ be the piecewise linear function
\begin{equation*}
\funcL(h) = \begin{cases}
\frac{h_0-h}{h_0-1} &,h\leq h_0 \\
0 &,h \geq h_0.
\end{cases}
\end{equation*}
\label{def:L}
\end{newdef}
We call a packing \emph{saturated} if no unit sphere can be added to the packing.
\begin{lem}[Lemma 6.95]
For any saturated packing $V$ and any $\mathbf{u}_0 \in V$,
\begin{equation}
\sum_{\mathbf{u}_1 \in V : \h\left(\mathbf{u}_0,\mathbf{u}_1\right)\leq h_0} \funcL\left(\h\left\lbrace\mathbf{u}_0,\mathbf{u}_1\right\rbrace\right)\leq 12.
\label{eq:lem6.95}
\end{equation}
\label{lem:6.95}
\end{lem}
We will not discuss the proof of this lemma since it is irrelevant to calculate the constant $c'$ explained in Section~\ref{sec:Introduction}. It is a computer proof anyways.
\begin{lem}[Lemma 6.97 ]
Inequality~(\ref{eq:lem6.95}) implies that for every saturated packing $V$, there exists a negligible FCC-compatible function $\G\colon V \to \mathbb{R}$.
\label{lem:6.97}
\end{lem}
Lemmas~\ref{lem:6.13},~\ref{lem:6.95}~and~\ref{lem:6.97} together imply the Kepler conjecture. As mentioned above, we are interested in replacing the constant $c$ in Lemma~\ref{lem:6.13} by a constant $c'$ that is independent of the packing $V$. In the proof of Lemma~\ref{lem:6.13} in \cite{Blueprint} it is shown that
\begin{equation}
\density(V,\mathbf{0},r)\leq \frac{\pi}{\sqrt{18}}\left(1+\frac{3}{r}\right)^3 + c_1\frac{\left(r+1\right)^2}{r^34\sqrt{2}},
\label{eq:lem:6.13}
\end{equation}
where $c_1$ is the constant from Definition~\ref{def:FCCnegl} for the FCC-compatible negligible function $\G\colon V\to \mathbb{R}$ that exists by assumption. So, what we need to show in order to find $c'$ is the following: There exists a constant $d$ independent of the function $\G$ and therefore of the packing $V$ such that we can set $c_1=d$ in Definition~\ref{def:FCCnegl} and Lemma~\ref{lem:6.97} still holds.

In the following, we derive such a constant $d$, reproduce the proof of Lemma~\ref{lem:6.97}, and show that the constructed FCC-compatible function also fulfills the stronger definition of negligible for $c_1=d$.

\section{Proof of Lemma~\ref{lem:6.97}}
For the following function, we will show that it is FCC-compatible and negligible:
\begin{equation}
\G(\mathbf{u})=-\vol \left(\voronoi\left(V,\mathbf{u}\right)\right)+8m_1-
\sum_{\mathbf{v}\in V\setminus\left\lbrace\mathbf{u}\right\rbrace} 8m_2 \funcL\left(\h\left(\left[\mathbf{u};\mathbf{v}\right]\right)\right)
\end{equation}
for constants $m_1$ and $m_2$ defined later and the function $\funcL$ as described in Definition~\ref{def:L}.

To prove that this function is FCC-compatible and negligible, another lemma is used. For this lemma we need more definitions.

In the following, if we talk about cells, we mean so called \emph{Marchal cells}. We refer the reader to Definition 6.51 in \cite{Blueprint} since we will not need most of the definition for the calculations made here. Cells are defined by four points in $V$ and a number $0\leq k \leq 4$ and are denoted by $\cell(\underline{\mathbf{u}},k)$ ("$k$-cell") for $\underline{\mathbf{u}}$ is a list of four points with some extra property to be explained later. 4-cells are always tetrahedra. The cells form a partition of $\mathbb{R}^3$. For a cell $X \neq \emptyset$, let $V(X) = X \cap V$ for a saturated packing $V$ (Definition 6.62 and Lemma 6.63 in \cite{Blueprint}).

The following three definitions build upon each other.
\begin{newdef}[Definition 3.7]
A set $C$ is \emph{r-radial} at center $\mathbf{v}$ if the two conditions $C \subset \B(\mathbf{v},r)$ and $\mathbf{v}+ \mathbf{u} \in C$ imply $\mathbf{v}+t\mathbf{u} \in C$ for all $t$ satisfying $0 < \Vert\mathbf{u}\Vert t < r$. A set $C$ is \emph{eventually radial} at center $\mathbf{v}$ if $C \cap \B(\mathbf{v},r)$ is r-radial at center $\mathbf{v}$ for some $r>0$.
\end{newdef}
\begin{newdef}[Definition 3.11 ]
When $C$ is measurable and eventually radial at center $\mathbf{v}$, define the \emph{solid angle} of $C$ at $\mathbf{v}$ to be 
\begin{equation*}
\sol(C,\mathbf{v}) = 3\frac{\vol\left(C \cap \B\left(\mathbf{v},r\right)\right)}{r^3},
\end{equation*}
where $r$ is as in the definition of eventually radial. By Lemma 3.10 in \cite{Blueprint}, this yields the same value when replacing r by $r'$ for any $0 \le r'\le r$. 
 \end{newdef}
\begin{newdef}[Definition 6.66]
Define the \emph{total solid angle} of a cell $X$ to be
\begin{equation*}
\tsol(X)=\sum_{\mathbf{v}\in V(X)} \sol(X,\mathbf{v}).
\end{equation*}
\end{newdef}
We will use this definition later. The following seven definitions build upon each other.
\begin{newdef}[Definition 5.43]
Recall that a set $A \subset \mathbb{R}^n$ is \emph{affine} if for every $\mathbf{v}, \mathbf{w} \in A$ and every $t \in \mathbb{R}$,
\begin{equation*}
t\mathbf{v} + \left( 1-t\right)\mathbf{w} \in A.
\end{equation*}
Recall that the \emph{affine hull} of $P \subset \mathbb{R}^n$ (denoted $\aff(P)$) is the smallest affine set containing $P$. The \emph{affine dimension} of $P$ (written $\dimaff (P)$ is $\card (S)-1$, where $S$ is a set of smallest cardinality such that $P \subset \aff(S)$. In particular, the affine dimension of the empty set is $-1$. [...]
\end{newdef}
\begin{newdef}[Definition 6.18]
Let $V$ be a saturated packing. When $k=0,1,2,3$, let $\underline{V}(k)$ be the set of lists $\underline{\mathbf{u}}=\left[\mathbf{u}_0; \dots; \mathbf{u}_k\right]$ of length $k+1$ with $\mathbf{u}_i \in V$ such that
\begin{equation*}
\dimaff \left(\Omega\left(V,\left[\mathbf{u}_0; \dots; \mathbf{u}_j\right]\right)\right)= 3-j
\end{equation*}
for all $0<j\leq k$. Set $\underline{V}(k)=\emptyset$ for $k >3$.
\label{def:6.18}
\end{newdef}
\begin{newdef}[Definition 6.24]
Let $V$ be a saturated packing and let $\underline{\mathbf{u}}=\left[\mathbf{u}_0; \dots ;\mathbf{u}_k\right] \in \underline{V}(k)$ for some $k$. Define points $\omega_j = \omega_j(V, \underline{\mathbf{u}}) \in \mathbb{R}^3$ by recursion over $j \le k$
\begin{align*}
\omega_0 &= \mathbf{u}_0\\
\omega_{j+1} &= \text{ the closest point to $\omega_j$ on $\voronoi\left(V, \left[\mathbf{u}_0; \dots ;\mathbf{u}_{j+1}\right]\right)$},
\end{align*}
Set $w(V,\underline{\mathbf{u}}) =\omega_k(V,\underline{\mathbf{u}})$, when $\underline{\mathbf{u}} \in \underline{V}(k)$. The set $V$ is generally fixed and is dropped from the notation. 
\end{newdef}
Let $\conv\lbrace \mathbf{v}_0,\dots, \mathbf{v_n}\rbrace$ denote the convex hull of the point set $\lbrace \mathbf{v}_0,\dots, \mathbf{v_n}\rbrace$.
\begin{newdef}[Definition 6.51]
Let $V$ be a saturated packing. Let $\underline{\mathbf{u}}=\left[\mathbf{u}_0; \dots ; \mathbf{u}_3\right] \in \underline{V}(3)$. Define $\xi(\underline{\mathbf{u}})$ as follows. If $\sqrt{2} \le \h\left(\left[\mathbf{u}_0;\dots ;\mathbf{u}_2\right]\right)$, then let $\xi(\underline{\mathbf{u}})=\omega\left(\left[\mathbf{u}_0;\dots ;\mathbf{u}_2\right]\right)$. If $\h\left(\left[\mathbf{u}_0;\dots ;\mathbf{u}_2\right]\right) < \sqrt{2} \le \h(\underline{\mathbf{u}})$, define $\xi(\underline{\mathbf{u}})$ to be the unique point in
\begin{equation*}
\conv \left\lbrace \omega\left(\left[\mathbf{u}_0;\dots ;\mathbf{u}_2\right]\right), \omega\left(\underline{\mathbf{u}}\right)\right\rbrace
\end{equation*}
at distance $\sqrt{2}$ from $\mathbf{u}_0$. [...]
\end{newdef}
\begin{newdef}[Definition 2.66]
When $\mathbf{v}_0 \neq \mathbf{v}_1$, write $\dih_V\left(\left\lbrace\mathbf{v}_0, \mathbf{v}_1\right\rbrace,\left\lbrace\mathbf{v}_2, \mathbf{v}_3\right\rbrace\right)$ for the angle $\gamma \in \left[0,\pi\right]$ formed by
\begin{align*}
\overline{\mathbf{w}}_2 &= (\mathbf{w}_1 \cdot \mathbf{w}_1) \mathbf{w}_2 - (\mathbf{w}_1 \cdot \mathbf{w}_2)\mathbf{w}_1 \text{ and}\\
\overline{\mathbf{w}}_3 &= (\mathbf{w}_1 \cdot \mathbf{w}_1) \mathbf{w}_3 - (\mathbf{w}_1 \cdot \mathbf{w}_3)\mathbf{w}_1,
\end{align*}
where $\mathbf{w}_i = \mathbf{v}_i - \mathbf{v}_0$.
\end{newdef}
\begin{newdef}[Definition 6.67]
Let $E(X)$ be the set of \emph{extremal edges} of the $k$-cell $X$ in a saturated packing $V$. More precisely, let
\begin{equation*}
E(X) = \left\lbrace\left\lbrace\mathbf{u}_i, \mathbf{u}_j\right\rbrace : \mathbf{u}_i \neq \mathbf{u}_j \in V(X) \right\rbrace.
\end{equation*}
In particular, $E(X)$ is empty for $0$ and $1$-cells and contains $\binom{k}{2}$ pairs when $2\le k \le 4$.
\label{def:6.67}
\end{newdef}
\begin{newdef}[Definition 6.68]
Let $V$ be a saturated packing. Let $X$ be a $k$-cell, where $2\le k \le 4$. Let $\varepsilon \in E(X)$. We define the dihedral angle $\dih(X, \varepsilon)$ of $X$ along $\varepsilon$ as follows. Explicitly, if $X$ is a null set, then set $\dih(X,\varepsilon)=0$. Otherwise, choose $\underline{\mathbf{u}}=[\mathbf{u}_0;\mathbf{u}_1;\mathbf{u}_2;\mathbf{u}_3] \in \underline{V}(3)$ such that $X=cell(\underline{\mathbf{u}},k)$ and $\varepsilon = \lbrace \mathbf{u}_0, \mathbf{u}_1\rbrace$. Set $\dih(X,\varepsilon)=\dih_V\left(\left\lbrace \mathbf{u}_0, \mathbf{u}_1\right\rbrace, \left\lbrace\mathbf{v}, \mathbf{w}\right\rbrace\right)$, where
\begin{align*}
\left\lbrace\mathbf{v},\mathbf{w}\right\rbrace =
\begin{cases}
\left\lbrace\xi(\underline{\mathbf{u}}), \omega(\underline{\mathbf{u}})\right\rbrace & ,k=2\\
\left\lbrace\mathbf{u}_2,\xi(\underline{\mathbf{u}})\right\rbrace & ,k=3\\
\left\lbrace\mathbf{u}_2, \mathbf{u}_3\right\rbrace &,k=4.
\end{cases}
\end{align*}
This is independent of the choice of $\underline{\mathbf{u}}$ defining $X$.
\end{newdef}
Now, we need just two more definitions before we can state the lemma mentioned above. The second one builds upon all the previously given definitions.
\begin{newdef}[Definition 6.70]
Define the following constants [...]:
\begin{align*}
\sol_0 &= 3\arccos\left(\frac{1}{3}\right)-\pi\\
\tau_0 &= 4\pi - 20\sol_0\\
m_1 &= \sol_0 2 \frac{\sqrt{2}}{\tau_0} \approx 1.012\\
m_2 & = \left(6\sol_0 - \pi\right)\frac{\sqrt{2}}{6\tau_0} \approx 0.0254\\
\end{align*}
[...]
\end{newdef}
\begin{newdef}[Definition 6.79]
For any cell $X$ of a saturated packing, define the function $\gamma(X,\ast)$ on $\left\lbrace f \colon \mathbb{R} \to \mathbb{R}\right\rbrace$ by
\begin{equation*}
\gamma(X,f) = \vol\left(X\right) - \left(\frac{2m_1}{\pi}\right)\tsol\left(X\right) + \left(\frac{8m_2}{\pi}\right)\sum_{\varepsilon \in E(X)} \dih\left(X,\varepsilon\right)f\left(h\left(\varepsilon\right)\right).
\end{equation*}
\end{newdef}
Now, we can state the lemma.
\begin{lem}[Lemma 6.86]
Let $f$ be any bounded, compactly supported function. Set
\begin{equation*}
\G\left(\mathbf{u}_0,f\right)=-\vol\left(\voronoi\left(V,\mathbf{u}_0\right)\right)+8m_1 - \sum_{\mathbf{u}\in V \setminus \left\lbrace \mathbf{u}_0\right\rbrace} 8m_2f\left(\h\left(\left[\mathbf{u}_0;\mathbf{u}\right]\right)\right).
\end{equation*}
If
\begin{equation*}
\sum_{\mathbf{v}\in V \setminus \left\lbrace\mathbf{u}\right\rbrace} f\left(\h\left(\left[\mathbf{u};\mathbf{v}\right]\right)\right) \le 12,
\end{equation*}
then $\G(\ast,f)$ is FCC-compatible. Moreover, if there exists a constant $c_0$ such that for all $r \ge 1$
\begin{equation*}
\sum_{X \subset \B\left(\mathbf{0},r\right)} \gamma\left(X,f\right) \ge c_0 r^2,
\end{equation*}
then $\G(\ast,f)$ is negligible. \footnote{In the sequel, by slight abuse of notation $\sum_{X \subset \B(\mathbf{0},r)}$ refers to the summation only over all Marchal cells $X \subset \B(\mathbf{0},r)$.}
\label{lem:6.86}
\end{lem}
This Lemma almost implies that there is a FCC-compatible negligible function. Using this Lemma, it is sufficient to show that Lemma~\ref{lem:6.95} implies the existence of suitable $\G$ and $f$ in Lemma~\ref{lem:6.86} to prove Lemma~\ref{lem:6.97}. We will turn to this later. Now, we focus on the proof of Lemma~\ref{lem:6.86}. Since we are only interested in the constant in the definition of negligible, we will only outline this part of the proof here. The idea of the proof of negligibility is as follows. If one can show that
\begin{align}
-\sum_{\mathbf{u} \in V(\mathbf{0},r)}\G\left(\mathbf{u},f\right) &\geq \sum_{X \subset \B(\mathbf{0},r)}\gamma\left(X,f\right) + c_2 r^2\label{eq:GandGamma}\\
\intertext{for some constant $c_2$. Then by the assumption in Lemma~\ref{lem:6.86}}\nonumber\\
-\sum_{\mathbf{u} \in V(\mathbf{0},r)}\G\left(\mathbf{u},f\right) &\geq \left(c_0 +c_2\right)r^2,\label{eq:negligible}
\end{align}
which directly implies, that $\G(\ast,f)$ is negligible for $c_1 = -(c_0 + c_2)$ in Definition~\ref{def:FCCnegl}. Since we are only interested in the independence of $c_1$ of the packing $V$, we will now focus on calculating the constant $c_2$ and thereby showing that $c_2$ is independent of the packing. We will get an estimate on $c_0$ later satisfying the conditions in Lemma~\ref{lem:6.86}.

\subsection{Calculation of $c_2$}
In this section we go through the proof of Lemma~\ref{lem:6.86} and show that 
\begin{equation*}
-\sum_{\mathbf{u} \in V(\mathbf{0},r)}\G\left(\mathbf{u},f\right) \geq \sum_{X \subset \B(\mathbf{0},r)}\gamma\left(X,f\right) + c_2 r^2
\end{equation*}
for some constant $c_2$ independent of the packing $V$. 

Observe that the equalities 
\begin{align*}
&-\sum_{\mathbf{u} \in V(\mathbf{0},r)}\G\left(\mathbf{u},f\right) = \\
&\sum_{\mathbf{u} \in V(\mathbf{0},r)}\vol\left(\voronoi\left(V,\mathbf{u}\right)\right)-\sum_{\mathbf{u} \in V(\mathbf{0},r)}8m_1 + \sum_{\mathbf{u} \in V(\mathbf{0},r)}\sum_{\mathbf{v}\in V \setminus \left\lbrace \mathbf{u}\right\rbrace} 8m_2f\left(\h\left(\left[\mathbf{u};\mathbf{v}\right]\right)\right)
\end{align*}
and
\begin{align*}
&\sum_{X \subset \B(\mathbf{0},r)}\gamma\left(X,f\right) = \\
&\sum_{X \subset \B(\mathbf{0},r)}\vol\left(X\right) - \sum_{X \subset \B(\mathbf{0},r)}\left(\frac{2m_1}{\pi}\right)\tsol\left(X\right) + \sum_{X \subset \B(\mathbf{0},r)}\left(\frac{8m_2}{\pi}\right)\sum_{\varepsilon \in E(X)} \dih\left(X,\varepsilon\right)f\left(h\left(\varepsilon\right)\right)
\end{align*}
have each three summands and we will relate them in this order. So, we start by showing that
\begin{align}
\sum_{\mathbf{u} \in V(\mathbf{0},r)}\vol\left(\voronoi\left(V,\mathbf{u}\right)\right) \geq \sum_{X \subset \B(\mathbf{0},r)}\vol\left(X\right) - \frac{56}{3}\pi r^2.\label{eq:GandGammFirst}
\end{align}
By Lemma 6.7 in \cite{Blueprint}, $\voronoi(V,\mathbf{v}) \subset \B(\mathbf{v},2)$, so we have
\begin{align*}
\sum_{\mathbf{u} \in V(\mathbf{0},r)}\vol\left(\voronoi\left(V,\mathbf{u}\right)\right) & \geq \vol\left(\B\left(\mathbf{0}, r-2\right)\right)\\
&=\vol\left(\B\left(\mathbf{0},r\right)\right) - \vol\left(\B\left(\mathbf{0},r\right)\setminus\B\left(\mathbf{0},r-2\right)\right)\\
&\geq \sum_{X \subset \B\left(\mathbf{0},r\right)}\vol\left(X\right) - \left(\frac{4}{3}\pi r^3 - \frac{4}{3}\pi\left(r-2\right)^3\right),\\
\intertext{since Marchal cells form a partition of space as mentioned above. So we get }
\sum_{\mathbf{u} \in V(\mathbf{0},r)}\vol\left(\voronoi\left(V,\mathbf{u}\right)\right) & \geq \sum_{X \subset \B\left(\mathbf{0},r\right)}\vol\left(X\right) - 8 \pi r^2 - \frac{32}{3}\pi.\\
\intertext{Since negligible is only defined for $r\geq 1$, we can assume $r\geq 1$ here. So we have}\\
\sum_{\mathbf{u} \in V(\mathbf{0},r)}\vol\left(\voronoi\left(V,\mathbf{u}\right)\right) &\geq \sum_{X \subset \B(\mathbf{0},r)}\vol\left(X\right) - \frac{56}{3}\pi r^2,
\end{align*}
which proves (\ref{eq:GandGammFirst}).

Next, we will show that
\begin{align}
-\sum_{\mathbf{u}\in V(\mathbf{0},r)} 8 m_1 \ge -\left(\frac{2m_1}{\pi}\right)\sum_{X \subset \B(\mathbf{0},r)}\tsol(X)-m_1 \cdot 2240 r^2.\label{eq:GandGammaSecond}
\end{align}
First, we need to estimate the diameter  of a Marchal cell. For a cell $X=\cell(\underline{\mathbf{u}},k)$ for $\underline{\mathbf{u}}=[\mathbf{u}_0;\dots] \in \underline{V}(3)$ it holds that
\begin{align*}
X = \cell(\underline{\mathbf{u}},k) &\subset \voronoi(V,\mathbf{u}_0) \cup \dots \cup \voronoi(V,\mathbf{u}_{k-1})\\
\intertext{(see the proof of Lemma 6.63 in \cite{Blueprint}). By Lemma 6.7, $\voronoi(V,\mathbf{v}) \subset \B(\mathbf{v},2)$, so we get}
X &\subset \B(\mathbf{u}_0,2) \cup \dots \cup \B(\mathbf{u}_{k-1},2).
\end{align*}
By Definition~\ref{def:6.18}, the Voronoi cells $\voronoi(V,\mathbf{u}_0), \dots, \voronoi(V,\mathbf{u}_{k-1})$ share at least one point and so the spheres $\B(\mathbf{u}_0,2), \dots, \B(\mathbf{u}_{k-1},2)$ do. Therefore, we can conclude that the diameter of a Marchal cell is at most 4.

Now,
\begin{align*}
\sum_{X\subset \B(\mathbf{0},r)}\tsol(X) &= \sum_{X\subset \B(\mathbf{0},r)}\sum_{\mathbf{v}\in V(X)} \sol(X,\mathbf{v})\\
\intertext{by definition. Since the diameter of each cell is at most 4, we can rewrite this as}
\sum_{X\subset \B(\mathbf{0},r)}\tsol(X)&= \sum_{X \subset \B(\mathbf{0},r+4)}\sum_{\mathbf{u}\in V(X)}\sol(X,\mathbf{u}) - \sum_{X \subset \B(\mathbf{0},r+4):X \nsubseteq
\B(\mathbf{0},r)}\sum_{\mathbf{u}\in V(X)}\sol(X,\mathbf{u})\\
&\ge \sum_{\mathbf{v}\in V(\mathbf{0},r)}\sum_{X:\mathbf{v}\in V(X)}\sol\left(X,\mathbf{v}\right) - \sum_{\mathbf{u}\in V(\mathbf{0},r+4)\setminus V(\mathbf{0},r-4)}\sum_{X:\mathbf{u}\in V(X)} \sol(X,\mathbf{u}),\\
\intertext{since $V(X) = X \cap V$ and a cell with diameter at most 4 that is not completely contained in a sphere of radius $r$ cannot intersect the sphere with the same center of radius $r-4$. Next, we use that the solid angles around one point sum up to $4\pi$. Furthermore, to estimate the number of points $\mathbf{u}\in V(\mathbf{0},r+4)\setminus V(\mathbf{0},r-4)$, we use the following volume argument. Each $\mathbf{u} \in V(\mathbf{0},r+4)\setminus V(\mathbf{0},r-4)$ is the center of a packed unit sphere and therefore the volume of those unit spheres sum up to less than the volume of $\B(\mathbf{0},r+5)\setminus \B(\mathbf{0},r-5)$. Vice versa dividing $\vol(\B(\mathbf{0},r+5)\setminus \B(\mathbf{0},r-5))$ by the volume of a unit sphere yields an upper bound on the number of points, so we get}
\sum_{X\subset \B(\mathbf{0},r)}\tsol(X) &\geq \sum_{\mathbf{v}\in V(\mathbf{0},r)} 4\pi - \left(\frac{\frac{4}{3}\pi\left(r+5\right)^3-\frac{4}{3}\pi\left(r-5\right)^3}{\frac{4}{3}\pi}\right)4\pi\\
&=  \sum_{\mathbf{v}\in V(\mathbf{0},r)} 4\pi - (30r^2 + 250)4\pi\\
&\ge  \sum_{\mathbf{v}\in V(\mathbf{0},r)} 4\pi - 1120\pi r^2,\\
\intertext{since $r\geq 1$.}
\end{align*}
By rearranging, we obtain
\begin{equation*}
-\left(\frac{2m_1}{\pi}\right)\sum_{X\subset\B(\mathbf{0},r)}\tsol(X) - 2240m_1 r^2 \le -\sum_{\mathbf{v}\in V(\mathbf{0},r)} 8m_1,
\end{equation*}
proving (\ref{eq:GandGammaSecond}).

Finally, we show that
\begin{equation}
\sum_{\mathbf{u}\in V(\mathbf{0},r)}\sum_{\mathbf{v}\in V \setminus \left\lbrace \mathbf{u}\right\rbrace} 8m_2f\left(\h\left(\left[\mathbf{u};\mathbf{v}\right]\right)\right) \ge \frac{8m_2}{\pi}\sum_{X\in\B(\mathbf{0},r)}\sum_{\varepsilon\in E(X)} \dih(X,\varepsilon)f\left(\h\left(\varepsilon\right)\right).\label{eq:GandGammThird}
\end{equation}
For any $\varepsilon \in E(X)$ it holds that $\varepsilon \subset V(X) \subset X$ by Definition~\ref{def:6.67}, and Lemma~6.63 in \cite{Blueprint}. This implies
\begin{align*}
\sum_{X\subset\B(\mathbf{0},r)}\sum_{\varepsilon\in E(X)} \dih(X,\varepsilon)f\left(\h\left(\varepsilon\right)\right) &\le \sum_{\varepsilon=\lbrace\mathbf{u},\mathbf{v}\rbrace \subset\B(\mathbf{0},r)}\sum_{X:\varepsilon\in E(X)} \dih\left(X,\varepsilon\right)f\left(\h\left(\varepsilon\right)\right)\\
&= \sum_{\varepsilon \subset \B(\mathbf{0},r)}2\pi f\left(\h\left(\varepsilon\right)\right),\\
\intertext{since the dihedral angle around an edge sums up to $2\pi$.\footnotemark When summing over ordered pairs, each edge appears twice, so we have to divide by two and get}
\sum_{X\subset\B(\mathbf{0},r)}\sum_{\varepsilon\in E(X)} \dih(X,\varepsilon)f\left(\h\left(\varepsilon\right)\right)
&\le \sum_{\mathbf{u}\in V(\mathbf{0},r)}\sum_{\mathbf{v}\in V(\mathbf{0},r)\setminus \lbrace \mathbf{u}\rbrace} \pi f\left(\h\left(\mathbf{u},\mathbf{v}\right)\right).
\end{align*}
\footnotetext{Here, we follow the argumentation of Hales \cite{Blueprint}. The author noticed that in Lemma 6.69 in \cite{Blueprint}, it is required that $\h(\varepsilon)<\sqrt{2}$ to have that the dihedral angles sum up to $2\pi$. Since Lemma~\ref{lem:6.86} is only used for $f=\funcL$ and $\funcL(\h(\varepsilon))=0$ for $\h(\varepsilon)\ge\sqrt{2}$, this does not cause any problems.}
So, for the last summand we get by rearrangement 
\begin{equation*}
\frac{8m_2}{\pi} \sum_{X\in\B(\mathbf{0},r)}\sum_{\varepsilon\in E(X)} \dih(X,\varepsilon)f\left(\h\left(\varepsilon\right)\right) \le \sum_{\mathbf{u}\in V(\mathbf{0},r)}\sum_{\mathbf{v}\in V(\mathbf{0},r)\setminus \lbrace \mathbf{u}\rbrace} 8m_2 f\left(\h\left(\mathbf{u},\mathbf{v}\right)\right),
\end{equation*}
showing (\ref{eq:GandGammThird}).

Combining (\ref{eq:GandGammFirst}), (\ref{eq:GandGammaSecond}), and (\ref{eq:GandGammThird}), we have
\begin{align*}
&-\G\left(\mathbf{u},f\right) \\
&=\sum_{\mathbf{u} \in V(\mathbf{0},r)}\vol\left(\voronoi\left(V,\mathbf{u}\right)\right)-\sum_{\mathbf{u} \in V(\mathbf{0},r)}8m_1 + \sum_{\mathbf{u} \in V(\mathbf{0},r)}\sum_{\mathbf{v}\in V \setminus \left\lbrace \mathbf{u}\right\rbrace} 8m_2f\left(\h\left(\left[\mathbf{u};\mathbf{v}\right]\right)\right) \\
&\ge\sum_{X \subset \B(\mathbf{0},r)}\vol\left(X\right) - \frac{56}{3}\pi r^2 -\left(\frac{2m_1}{\pi}\right)\sum_{X \subset \B(\mathbf{0},r)}\tsol(X)-m_1 \cdot 2240 r^2 + \\
&\frac{8m_2}{\pi}\sum_{X\in\B(\mathbf{0},r)}\sum_{\varepsilon\in E(X)} \dih(X,\varepsilon)f\left(\h\left(\varepsilon\right)\right)\\
&=\sum_{X \subset\B(\mathbf{0},r)}\gamma(X,f) - \left(\frac{56}{3}+m_1\cdot 2240\right)r^2,
\intertext{i.e. (\ref{eq:GandGamma}) holds with $c_2 = -\frac{56}{3}-m_1\cdot 2240$.}
\end{align*}

As mentioned above, we only need to show that the conditions in Lemma~\ref{lem:6.86} hold for some $f$ to prove Lemma~\ref{lem:6.97}. The first condition holds for $f=\funcL$ by Lemma~\ref{lem:6.95}. It remains to show that the second condition also holds for $f=\funcL$, namely $\sum_{X \subset \B\left(\mathbf{0},r\right)} \gamma\left(X,\funcL\right) \ge c_0 r^2$. Then, we will have showed that the function $\G(\ast,\funcL)$ is FCC-compatible and negligible for $c_1=-(c_0 + c_2)$ in Definition~\ref{def:FCCnegl}. For $c_2$ we already showed that it is independent of the packing, now, we will do so for $c_0$, i.e. prove that $\sum_{X \subset \B\left(\mathbf{0},r\right)} \gamma\left(X,\funcL\right) \ge c_0 r^2$ holds for a $c_0$ independent of the packing.

\subsection{Calculation of $c_0$}
In this section, we will show that
\begin{equation}
\sum_{X \subset \B\left(\mathbf{0},r\right)} \gamma\left(X,f\right) \ge c_0 r^2\label{eq:Gamma}
\end{equation}
holds for $f=\funcL$ for a constant $c_0$ independent of the packing that defines the cells $X$. We will need more definitions for this proof and we will try to give them as closely as possible to the point where they are used.
\begin{newdef}[Definition 6.70 and 6.88]
Set
\begin{align*} 
h_+ &= 1.3254.\\
\intertext{Let $\M \colon \mathbb{R} \to \mathbb{R}$ be the piecewise polynomial function}
\M (h)&= \begin{cases}
\frac{\sqrt{2}-h}{\sqrt{2}-1}\frac{h_+-h}{h_+-1}\frac{17h-9h^2-3}{5} &,h \le \sqrt{2}\\
0 &,h>\sqrt{2}.
\end{cases}\\
\intertext{Let $h_-\approx1.23175$ be the unique root of the quadratic polynomial $\M(h)-\funcL(h)$ that lies in the interval $[1.231,1.232]$.}
\end{align*}
\end{newdef}
\begin{newdef}[Definition 6.89]
A \emph{critical edge} $\varepsilon$ of a saturated packing $V$ is an unordered pair that appears as an element of $E(X)$ for some $k$-cell $X$ of the packing $V$ such that $h(\varepsilon) \in [h_-,h_+]$. Let $\EC(X)$ be the set of critical edges that belong to $E(X)$. If $X$ is any cell such that $\EC(X)$ is not empty, let the \emph{weight} $\wt(X)$ of $X$ be $1/\operatorname{card}(\EC(X))$.
\label{def:critEdge}
\end{newdef}
Now, we can start with estimating $\sum_{X \subset \B(\mathbf{0},r)} \gamma(X,\funcL)$.
\begin{align}
\sum_{X \subset \B(\mathbf{0},r)} \gamma(X,\funcL) &= \sum_{X \subset \B(\mathbf{0},r) : \EC(X) \neq \emptyset} \gamma(X,\funcL) + \underbrace{\sum_{X \subset \B(\mathbf{0},r): \EC(X) = \emptyset} \gamma(X,\funcL)}_{\geq 0 \text{ by Lemma 6.92 in \cite{Blueprint}}}\nonumber\\
&\ge \sum_{X \subset \B(\mathbf{0},r) : \EC(X) \neq \emptyset} \gamma(X,\funcL)\nonumber\\
&=  \sum_{X \subset \B(\mathbf{0},r) : \EC(X) \neq \emptyset} \left(\gamma(X,\funcL) \underbrace{\sum_{\varepsilon \in \EC(X)}\wt(X)}_{=1 \text{ by Definition~\ref{def:critEdge}}}\right)\nonumber\\
&= \sum_{X \subset \B(\mathbf{0},r) : \EC(X) \neq \emptyset}\sum_{\varepsilon \in \EC(X)}\gamma(X,\funcL)\wt(X)\nonumber\\
&= \sum_{\varepsilon \subset \B(\mathbf{0},r)}\sum_{X:\varepsilon \in \EC(X)} \gamma(X,\funcL)\wt(X) - \underbrace{\sum_{\varepsilon \subset \B(\mathbf{0},r)}\sum_{X:\varepsilon \in \EC(X), X \nsubseteq \B(\mathbf{0},r)} \gamma(X,\funcL)\wt(X)}_{=:\alpha}\label{eq:gammaIntermediate}
\end{align}
Next, we will show that the term $\alpha$ is upper bounded by a constant~times~$r^2$.
\begin{align}
\alpha &= \sum_{\varepsilon \subset \B(\mathbf{0},r)}\sum_{X:\varepsilon \in \EC(X), X \nsubseteq \B(\mathbf{0},r)} \gamma(X,\funcL)\wt(X)\nonumber\\
&\le \sum_{X \subset \B(\mathbf{0},r+4) \setminus \B(\mathbf{0},r-4)} \sum_{\varepsilon \in \EC(X)} \gamma(X,\funcL)\wt(X),\label{eq:1}
\intertext{since the diameter of $X$ is at most 4 and $X$ intersects the boundary of $\B(\mathbf{0},r)$. By (\ref{eq:1}),}
\alpha &\le \sum_{X \subset \B(\mathbf{0},r+4) \setminus \B(\mathbf{0},r-4)}\left(\gamma(X,\funcL) \underbrace{\sum_{\varepsilon \in \EC(X)} \wt(X)}_{=1 \text{ by Definition~\ref{def:critEdge}}}\right).\label{eq:2}
\intertext{Next, we plug in the definition of $\gamma(X,L)$. By (\ref{eq:2}),}
\alpha &\le\sum_{X \subset \B(\mathbf{0},r+4) \setminus \B(\mathbf{0},r-4)} 
\left(\vol\left(X\right) - \underbrace{\left(\frac{2m_1}{\pi}\right)\tsol\left(X\right)}_{\ge 0 \text{ by definition}} + \left(\frac{8m_2}{\pi}\right)\sum_{\varepsilon \in E(X)} \dih\left(X,\varepsilon\right)\funcL\left(h\left(\varepsilon\right)\right)\right)\nonumber\\
&\le \sum_{X \subset \B(\mathbf{0},r+4) \setminus \B(\mathbf{0},r-4)} \vol\left(X\right)+ \sum_{X \in \B(\mathbf{0},r+4) \setminus \B(\mathbf{0},r-4)}\left(\frac{8m_2}{\pi}\right)\sum_{\varepsilon \in E(X)} \dih\left(X,\varepsilon\right)\funcL\left(h\left(\varepsilon\right)\right).\label{eq:3}
\intertext{Since Marchal cells are a partition of the space, the volume of the cells contained in a spherical shell sums up to at most the volume of the spherical shell. Furthermore, it holds for any $\varepsilon \in E(X)$ that $\varepsilon \subset V(X) \subset X$, so we get from (\ref{eq:3})}
\alpha &\le \frac{4}{3}\pi\left(r+4\right)^3-\frac{4}{3}\pi\left(r-4\right)^3 + \frac{8m_2}{\pi}\sum_{\varepsilon\subset \B(\mathbf{0},r+4) \setminus \B(\mathbf{0},r-4)}\sum_{X:\varepsilon\in E(X)}\dih\left(X,\varepsilon\right)\funcL\left(h\left(\varepsilon\right)\right).\label{eq:4}
\intertext{Since the dihedral angle for edges $\varepsilon$ with $\h(\varepsilon)<\sqrt{2}$ sums up to $2\pi$ (see Lemma~6.69 in \cite{Blueprint}) and for $\h(\varepsilon)\ge \sqrt{2}$ the factor $\funcL(\h(\varepsilon))=0$, we obtain from (\ref{eq:4})}
\alpha &\le 32\pi r^2 + \frac{512}{3} \pi + \sum_{\varepsilon\subset \B(\mathbf{0},r+4) \setminus \B(\mathbf{0},r-4)}16m_2\funcL\left(h\left(\varepsilon\right)\right).\label{eq:5}
\intertext{Now, we want to estimate the number of edges such that $\funcL(\h(\varepsilon)) > 0$. First, we give an upper bound on the number of points in $V$ inside the spherical shell as follows. By decreasing the inner radius and increasing the outer radius by 1, all unit spheres packed with centers in the original shell are completely contained in the enlarged shell. Then, we divide the volume of the enlarged shell by the volume of a unit sphere. Next, we want to give  an upper bound on the number of points in $V$ with distance at most 2.52 from a given point, since for longer edges $\funcL(\h(\varepsilon))=0$. We do this in a similar way as before by a volume argument. By multiplying these two values, we count each edge twice, so we have to divide by 2. In addition, the function $\funcL$ is upper bounded by $1.26/0.26$. So, by (\ref{eq:5})}
\alpha&\le 32\pi r^2 + \frac{512}{3} \pi +\underbrace{\frac{\frac{4}{3}\pi\left(r+5\right)^3-\frac{4}{3}\pi\left(r-5\right)^3}{\frac{4}{3}\pi}}_{\ge | V \cap \left(\B(\mathbf{0},r+4) \setminus \B(\mathbf{0},r-4)\right)|}\cdot\underbrace{\frac{\frac{4}{3}\pi \cdot 3.52^3}{\frac{4}{3}\pi}}_{\ge |V \cap \B(\mathbf{v},2.52)| \text{ for any }\mathbf{v}}\cdot \frac{1}{2}\cdot16m_2\cdot\frac{126}{26}\nonumber\\
&\le 32\pi r^2 + \frac{512}{3}\pi + \left(30r^2 + 250\right)\cdot 3.52^3\cdot 8 \cdot \frac{63}{13}\cdot 0.0255\nonumber\\
&\le 1394.1 r^2 + 11315.6\nonumber
\intertext{As mentioned before, we can assume $r\ge 1$ and obtain}
\alpha&\le 12710 r^2.\label{eq:alpha}
\end{align}
Now, we know from inequalities (\ref{eq:gammaIntermediate}) and (\ref{eq:alpha}) that
\begin{align}
\sum_{X \subset\B(\mathbf{0},r)} \gamma(X,L) \ge \sum_{\varepsilon \subset \B(\mathbf{0},r)}\sum_{X:\varepsilon \in \EC(X)} \gamma(X,L)\wt(X) - 12710^2.\label{eq:gammaIntermediate2}
\end{align}
To proceed with the estimation, we need two more definitions.
\begin{newdef}[Definition 6.90]
Set
\begin{equation*}
\beta_0(h) = 0.005\left(1-\frac{\left(h- h_0\right)^2}{\left(h_+ - h_0\right)^2}\right).
\end{equation*}
If $X$ is a 4-cell with exactly two critical edges and if those edges are opposite, then set
\begin{equation*}
\beta(\varepsilon,X) = \beta_0\left(\h\left(\varepsilon\right)\right) - \beta_0\left(\h\left(\varepsilon'\right)\right), \text{ where } \EC(X) = \lbrace\varepsilon,\varepsilon'\rbrace.
\end{equation*}
Otherwise, for all other edges in all other cells, set $\beta(\varepsilon,X)=0$.
\end{newdef}
\begin{newdef}[Definition 6.91]
Let $V$ be a saturated packing. Let $\varepsilon \in \EC(X)$ be a critical edge of a $k$-cell $X$ of $V$ for some $2\le k\le 4$. A \emph{cell cluster} is the set
\begin{equation*}
\CL(\varepsilon) = \left\lbrace X : \varepsilon \in \EC(X)\right\rbrace
\end{equation*}
of all cells around $\varepsilon$. Define
\begin{equation*}
\Gamma(\varepsilon) = \sum_{X \in \CL(\varepsilon)} \gamma(X,L)\wt(X) + \beta(\varepsilon,X).
\end{equation*}
\end{newdef}
Using these definitions, we can rewrite inequality (\ref{eq:gammaIntermediate2}) as follows.
\begin{align}
\sum_{X \subset\B(\mathbf{0},r)} \gamma(X,L) &\ge \sum_{\varepsilon \subset \B(\mathbf{0},r)}\sum_{X:\varepsilon \in \EC(X)} \gamma(X,L)\wt(X) - 12710^2\nonumber\\
&=\sum_{\varepsilon \subset \B(\mathbf{0},r)} \left(\Gamma(\varepsilon) - \sum_{X \in \CL(\varepsilon)}\beta(\varepsilon,X)\right) - 12710r^2\nonumber\\
&=\sum_{\varepsilon \subset \B(\mathbf{0},r)} \Gamma(\varepsilon) - \underbrace{\sum_{\varepsilon \subset \B(\mathbf{0},r)}\sum_{X \in \CL(\varepsilon)}\beta(\varepsilon,X)}_{=: \zeta}- 12710r^2\label{eq:gammaIntermediate3}
\end{align}
Next, we will upper bound the term $\zeta$.
\begin{align}
\zeta &= \sum_{\varepsilon \subset \B(\mathbf{0},r)}\sum_{X \in \CL(\varepsilon)}\beta(\varepsilon,X)\nonumber\\
&= \sum_{\varepsilon \subset \B(\mathbf{0},r)}\sum_{X \in \CL(\varepsilon): X \subset \B(\mathbf{0},r)}\beta(\varepsilon,X) + \sum_{\varepsilon \subset \B(\mathbf{0},r)}\sum_{X \in \CL(\varepsilon):X \nsubseteq \B(\mathbf{0},r)}\beta(\varepsilon,X)\nonumber\\
&= \sum_{X \subset \B(\mathbf{0},r)}\sum_{\varepsilon \in \EC(X)}\beta(\varepsilon,X) + \sum_{\varepsilon \subset \B(\mathbf{0},r)}\sum_{X \in \CL(\varepsilon):X \nsubseteq \B(\mathbf{0},r)}\beta(\varepsilon,X)\nonumber\\
&= 0+\sum_{\varepsilon \subset \B(\mathbf{0},r)}\sum_{\substack{X \in \CL(\varepsilon):X \nsubseteq \B(\mathbf{0},r)\\ \wedge X \text{ is 4-cell with}|\EC(X)|=2}}\beta(\varepsilon,X),\label{eq:6}
\intertext{since $\beta(\varepsilon,X) > 0$ only for 4-cells with exactly two critical edges $\varepsilon, \varepsilon'$ and $\beta(\varepsilon,X) + \beta(\varepsilon',X)=0$. Furthermore, $\beta(\varepsilon,X) \le 0.005$, so we get from (\ref{eq:6})}\nonumber\\
 \zeta &\le \sum_{\varepsilon \subset \B(\mathbf{0},r)}\sum_{\substack{X \in \CL(\varepsilon):X \nsubseteq \B(\mathbf{0},r)\\ \wedge X \text{ is 4-cell with }|\EC(X)|=2}} 0.005. \label{eq:7}
\intertext{Next, we want to exchange the inner and outer sum. Since the inner sum is zero for non-critical edges and the cells intersect the boundary of $\B(\mathbf{0},r)$, we get from (\ref{eq:7})}\nonumber\\
\zeta &\le\sum_{\substack{X \subset \B(\mathbf{0},r+4) \setminus \B(\mathbf{0},r-4) :\\ X \text{ is a 4-cell with }|\EC(X)|=2}}\sum_{\varepsilon\in \EC(X)} 0.005\nonumber\\
&\le\sum_{\substack{X \subset \B(\mathbf{0},r+4) \setminus \B(\mathbf{0},r-4) : \\X \text{ is a 4-cell}}} 0.01.\label{eq:8}
\intertext{Now, we will estimate the number of points in $V$ inside the sperical shell as before. A 4-cell is the convex hull of four points in $V$ with circumradius at most $\sqrt{2}$ (see Definition~6.51 in \cite{Blueprint}). Therefore, for a fixed point $\mathbf{v}\in V$ all points that can form a 4-cell with $\mathbf{v}$ must lie inside a ball of radius $2\sqrt{2}$ centered at $\mathbf{v}$. As before, we calculate an upper bound on the number of points in $V$ inside a ball of radius $2\sqrt{2}$. Then, the number of subsets of size three that can be formed with these points is an upper bound on the number of 4-cells a point in $V$ can be part of. By multiplying these numbers, we count each 4-cell 4 times, so we divide by 4. So, we get from (\ref{eq:8})}\nonumber\\
\zeta & \le \frac{\frac{4}{3}\pi(r+5)^3 - \frac{4}{3}\pi(r-5)^3}{\frac{4}{3}\pi} \cdot {\left\lfloor\frac{\frac{4}{3}\pi(2\sqrt{2}+1)^3}{\frac{4}{3}\pi}\right\rfloor \choose 3} \cdot \frac{1}{4}\cdot 0.01\nonumber\\
& = \frac{\left(30r^2 + 250\right)}{4}\cdot {\left\lfloor\left( 2\sqrt{2}+1\right)^3\right\rfloor \choose 3}\cdot 0.01\nonumber\\
&= \left(7.5 r^2 + 62.5\right)\cdot 27720\cdot 0.01\nonumber\\
&= 2079r^2 + 17325.\nonumber\\
\intertext{Again, we can assume $r\ge 1$, so}\nonumber
\zeta &\le 19404 r^2.
\end{align}
Now, we plug this into inequality~(\ref{eq:gammaIntermediate3}) and obtain
\begin{align*}
\sum_{X \subset\B(\mathbf{0},r)} \gamma(X,L) &\ge \sum_{\varepsilon \subset \B(\mathbf{0},r)} \Gamma(\varepsilon) - 19404r^2- 12710r^2\\
&= \underbrace{\sum_{\varepsilon \subset \B(\mathbf{0},r)} \Gamma(\varepsilon)}_{\ge 0 \text{ by Theorem~6.93 in \cite{Blueprint}}}-32114r^2\\
&\ge - 32114r^2,
\end{align*}
i.e. (\ref{eq:Gamma}) holds for $c_0= -32114$.

We showed that both conditions in Lemma~\ref{lem:6.86} hold and therefore the function $\G(\ast,\funcL)$ is FCC-compatible and negligible (see (\ref{eq:negligible}) and the explanation thereafter) for 
\begin{align*}
c_1=-(c_0 + c_2) &= \frac{56}{3}+m_1\cdot 2240 + 32114\\
 &\le  \frac{56}{3}+1.013\cdot 2240 + 32114\\
 &\le 34402.
\end{align*}
By (\ref{eq:lem:6.13}), the constant in Lemma~\ref{lem:6.13} only depends on constants and $c_1$. So it is independent of the packing since we showed that $c_1$ is. Therefore, we turn now to inequality (\ref{eq:lem:6.13}) and will later plug in $c_1$ as calculated above.
\begin{align*}
\density(V,\mathbf{0},r)&\leq \frac{\pi}{\sqrt{18}}\left(1+\frac{3}{r}\right)^3 + c_1\frac{\left(r+1\right)^2}{r^34\sqrt{2}}\\
&= \frac{\pi}{\sqrt{18}}\left(1+\frac{9}{r}+\frac{27}{r^2}+\frac{27}{r^3}\right)+c_1\frac{r^2+2r+1}{r^3 4\sqrt{2}}\\
&=\frac{\pi}{\sqrt{18}} + \frac{\pi}{\sqrt{18}}\left(\frac{9}{r}+\frac{27}{r^2}+\frac{27}{r^3}\right) + c_1\left(\frac{1}{r4\sqrt{2}}+\frac{2}{r^24\sqrt{2}}+\frac{1}{r^34\sqrt{2}}\right)
\intertext{Since $r\ge 1$, we have $\frac{1}{r^3} \le \frac{1}{r^2} \le \frac{1}{r}$ and get}
&\le \frac{\pi}{\sqrt{18}} + \frac{63\pi}{\sqrt{18}r}+\frac{c_1}{\sqrt{2}r}\\
&= \frac{\pi}{\sqrt{18}} + \frac{21\pi+c_1}{\sqrt{2}r}\\
&=\frac{\pi}{\sqrt{18}} + \frac{21\pi+34402}{\sqrt{2}}\cdot\frac{1}{r}.
\end{align*}

Summarizing, we showed that the constant in Lemma~\ref{lem:6.13} does not depend on the particular packing $V$ but only on the constant for the assumed existing FCC-compatible negligible function. Then, we showed that there is a FCC-compatible negligible function for which the definition for negligible holds for a constant independent of the packing. So, we can state the main result of this work as follows.
\begin{thm}
For a saturated packing $V$ and all $r\ge 1$ it holds that
\begin{equation*}
\density(V,\mathbf{0},r) \le \frac{\pi}{\sqrt{18}} + \frac{21\pi+34402}{\sqrt{2}}\cdot\frac{1}{r} \le \frac{\pi}{\sqrt{18}} + 24373\cdot\frac{1}{r}.
\end{equation*}
\end{thm}

\bibliographystyle{hplain}
	\bibliography{references2}
		
\end{document}